\documentclass{article}

\usepackage[a4paper, left=3.3cm,top=3.3cm,right=3.3cm,bottom=3.7cm]{geometry}
\usepackage{concmath}
\usepackage[T1]{fontenc}

\usepackage{amsmath,amsfonts,amssymb}
\usepackage{graphicx}
\usepackage{dsfont}
\usepackage{mathtools}
\usepackage{siunitx}
\usepackage{psfrag}
\usepackage{floatrow}

\usepackage{authblk}
\title{Photoacoustic image reconstruction via deep learning}

\author{Stephan Antholzer}

\affil{Department of Mathematics, University of Innsbruck\authorcr
Technikerstrasse 13, 6020 Innsbruck, Austria\authorcr
{\tt  stephan.antholzer@uibk.ac.at}}

\author{Johannes Schwab}

\affil{Department of Mathematics, University of Innsbruck\authorcr
Technikerstrasse 13, 6020 Innsbruck, Austria\authorcr
{\tt johannes.schwab@uibk.ac.at}}

\author{Robert Nuster}

\affil{Department of Physics, Universit\"at Graz\authorcr
Universitaetsplatz 5, Graz, Austria\authorcr
E-mail: {\tt  ro.nuster, guenther.paltauf@uni-graz.at }}

\author{Markus Haltmeier}

\affil{Department of Mathematics, University of Innsbruck\authorcr
Technikerstrasse 13, 6020 Innsbruck, Austria\authorcr
E-mail: {\tt markus.haltmeier@uibk.ac.at}}

\date{Januar 30, 2018}

\newcommand*{\R}{\mathds{R}}

\newcommand{\source}{p_0}
\newcommand{\wave}{\mathcal{W}}

\newcommand{\curve}{{\boldsymbol S}}
\newcommand{\rr}{\mathbf{r}}
\newcommand{\rrs}{\mathbf{s}}

\newcommand{\fbp}{\mathcal{B}}

\newcommand*{\Int}[4]{\int_{#1}^{#2}\!{#3}\,\mathrm{d}{#4}}

\newcommand*{\W}{\mathds{W}}
\newcommand*{\FBP}{\mathds{B}}
\newcommand*{\CNN}{\boldsymbol{\Phi}_{\mathrm{CNN}}}
\newcommand*{\UNET}{\boldsymbol{\Phi}_{\textnormal{U-Net}}}
\newcommand*{\SNET}{\boldsymbol{\Phi}_{\textnormal{S-Net}}}
\newcommand*{\NN}{\boldsymbol{\Phi}}
\newcommand*{\waved}{\mathds{A}}
\newcommand{\Dnum}{\mathds{D}}
\newcommand\set[1]{\left\{#1\right\}}

\newcommand*{\Ntrain}{N_{\TT}}
\newcommand*{\Nx}{N}
\newcommand*{\Nt}{M}
\newcommand*{\Ny}{d}
\newcommand*{\Xin}{\mathbf{Y}}
\newcommand*{\Xout}{\mathbf{X}}
\newcommand*{\noise}{\boldsymbol{\xi}}
\newcommand*{\TT}{\mathcal{T}}
\newcommand*{\WW}{\mathcal{W}}

\newcommand*{\MSE}{\ensuremath{\operatorname{MSE}}}

\renewcommand{\phi}{\varphi}
\renewcommand{\epsilon}{\varepsilon}

\DeclarePairedDelimiter{\abs}{\lvert}{\rvert}
\DeclarePairedDelimiter{\norm}{\lVert}{\rVert}

\DeclareMathOperator{\ReLU}{ReLU}

\usepackage{enumitem}
\setitemize{label=\scriptsize{$\blacksquare$}, wide}

\usepackage{soul}
\usepackage{xcolor}
\colorlet{lred}{red!40}
\colorlet{lgreen}{green!40}
\colorlet{lblue}{blue!40}

\begin{document}
\maketitle

\begin{abstract}
Applying standard algorithms to sparse data  problems in photoacoustic tomography (PAT) yields low-quality images containing severe under-sampling artifacts.  To some extent, these artifacts can be reduced by iterative image reconstruction algorithms which allow to include prior knowledge such as smoothness, total variation (TV) or sparsity constraints. These algorithms tend to be time consuming as the forward and adjoint problems have to be solved repeatedly. Further, iterative algorithms have additional drawbacks. For example, the reconstruction quality strongly depends on  a-priori  model assumptions about  the objects to be recovered, which are often not strictly satisfied in practical applications.
To overcome these issues, in this paper, we develop  direct and efficient reconstruction algorithms based on  deep learning. As opposed to iterative algorithms, we apply a convolutional neural network, whose  parameters are trained  before the reconstruction process  based on a set of training data. For actual image reconstruction, a single evaluation of the trained network yields the desired result.
Our presented numerical results (using two different network architectures)  demonstrate that the proposed deep learning approach
reconstructs images with a quality comparable to state of the art  iterative reconstruction methods.

\bigskip\noindent
\textbf{Keywords:}
Photoacoustic tomography, sparse data,  limited view problem, image reconstruction, deep learning, convolutional neural networks,  inverse problems.

\end{abstract}

\section{Introduction}
\label{sec:intro}

Deep learning is a rapidly emerging research topic, improving the performance of many image processing and computer vision systems.
Deep learning uses a rich class of learnable functions in the form of artificial neural networks, and contain free parameters that can be adjusted to the particular problem at hand. It is state of the art in many different domains and outperforms most comparable algorithms~\cite{goodfellow2016deep}{}.
However, only recently they have been used for image reconstruction, see
for example ~\cite{antholzer2017deep,chen2017lowdose,kelly2017deep,han2016deep,jin2017deep,wang2016perspective,zhang2016image}{}.
In image reconstruction with deep learning, a  convolutional neural network (CNN)
is designed to map the measured data to the desired reconstructed image.
The CNN contains free weights that can be adjusted   prior to the actual image reconstruction based on an appropriate  set of training images.  Actual image reconstruction with deep learning consists in a single  evaluation of the trained network. Instead of providing an explicit a-priori model, deep learning uses  training data and the network itself adjusts to the appropriate image reconstruction task and phantom class.

In this paper we develop  deep learning  based image reconstruction algorithms
for photoacoustic tomography (PAT), a hybrid imaging method for visualizing light
absorbing structures within optically scattering media. The image reconstruction task in PAT is to recover the internal
absorbing structures  from acoustic measurements made outside of the sample (see Figure~\ref{fig:PATPrinciple}).
We solve the PAT image reconstruction  problem by first applying the
filtered back-projection (FBP) algorithm to the measured data and subsequently applying  a trained CNN.
There  are plenty of existing CNN  architectures that can be combined with the FBP.
In  the present  paper we compare two different CNNs for that purpose.
The first one is (a slight variant  of) the U-Net developed  in \cite{ronneberger2015unet}
for image segmentation and winner of several machine learning competitions.
For comparison purpose, we also test a  self-designed  very simple CNN, named S-Net, that  consists of just   three convolutional layers.  Numerical results demonstrate
that both networks work well for PAT image reconstruction. It might be surprising  that  the  basic S-Net
already performs that well and yields a reconstruction quality comparable to the U-Net. The  design of  other simple network architectures (that can be evaluated faster than more complex ones)   even  outperforming  the U-Net for PAT is an interesting future challenge.

Using deep learning and in particular deep CNNs for image reconstruction in PAT has first  been proposed  in
our previous work~\cite{antholzer2017deep}{}. In particular, in that paper we proposed the combination of the the FBP with a trained network
for which the U-Net has been used.
Other learning approaches to PAT can be found in\cite{reiter2017machine,dreier2017operator,hauptmann2017model,schwab2018dalnet}{}.
Opposed to  \cite{antholzer2017deep}{}, in this paper we present numerical results using the S-Net and compare it to the U-Net.
Moreover, we use a different class  of test phantoms (for training and testing) that mimic a nonuniform, more realistic,   light distribution
within the samples, and different measurement setups including limited view.

\paragraph{Outline:}
The remainder of this paper is organized as follows. In Section~\ref{sec:pat} we give
a  brief introduction to PAT and formally  describe the  considered PAT image reconstruction problem. Deep learning  with an emphasis on its use
for image reconstruction is  presented  Section~\ref{sec:deep}. Our proposed
 deep learning approach for  PAT image reconstruction and the proposed network designs
 are   presented in  Section~\ref{sec:dpat}. Details of our numerical studies  and numerical results  are presented in Section~\ref{sec:results}. Thereby we also describe the network training and the considered training and evaluation data. A short summary and conclusions are presented in Section~\ref{sec:conclusion}.

\begin{psfrags}
\psfrag{A}{$\curve$}
\psfrag{B}{}
\begin{figure}[htb!]
    \center
    \includegraphics[width=\textwidth]{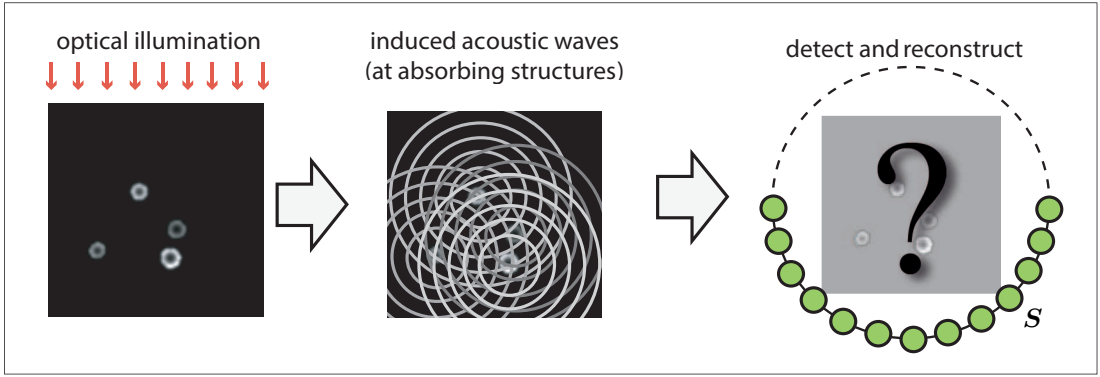}
    \caption{\textbf{Basic  principles of PAT.} Left: A sample object is illuminated by short optical  pulses. Middle: Optical energy is absorbed within the sample, causes nonuniform heating and induces a subsequent  acoustic pressure wave. Right: Acoustic sensors located outside of the sample
    capture the pressure signals, which are used to recover an image of the interior. In this paper we use deep learning and in particular deep CNNs
    for image reconstruction. Our approach allows a small number of sensor positions arranged on a possibly non-closed measurement curve $\curve$.} \label{fig:PATPrinciple}
\end{figure}
\end{psfrags}

\section{Methods}

\subsection{Photoacoustic tomography}
\label{sec:pat}

PAT is a non-invasive coupled-physics biomedical imaging technique offering high contrast and high spatial resolution~\cite{kruger1995photoacoustic,paltauf2007photacoustic}{}.
As illustrated in Figure\ref{fig:PATPrinciple},   a semi-transparent sample is illuminated
with short optical pulses which causes heating of the sample followed by expansion and the subsequent emission of an acoustic pressure wave.
Detectors outside of the sample measure the acoustic wave and the  measurements are then used to reconstruct the initial pressure, which provides information about the interior of the object.
We denote the initial pressure distribution  (which is the function to be reconstructed) by  $\source \colon \R^d \to \R$.
The cases  $d=2$ and $d=3$ for the spatial dimension are relevant in applications.
In order to simplify the presentation, in the following   we only consider the case $d=2$.
The 2D case arises, for example, when one uses integrating line detectors in PAT, see \cite{burgholzer2007temporal,paltauf2007photacoustic}{}.

In two spatial dimensions, the induced acoustic pressure $p\colon\R^2\times [0,\infty)\to \R$ in PAT satisfies the following initial value problem
\begin{equation} \label{eq:wave}
\left\{\begin{array}{ll}
    \partial^2_t p (\rr,t)  - \Delta_{\rr} p(\rr,t) = 0 & \text{ for } (\rr,t) \in \R^2 \times (0,\infty) \\
p(\rr,0) = \source(\rr)   & \text{ for } \rr \in \R^2 \\
\partial_t p(\rr,0) = 0 & \text{ for } \rr \in \R^2 \,.
\end{array}\right.
\end{equation}
Here $\rr\in \R^2$ is the  spatial  variable, $\Delta_{\rr}$ the spatial Laplacian, $t$ the (rescaled) time variable,
and $\partial_t$ the temporal derivative.  We assume that the sound speed $v_s$  is constant and that the physical time variable
$ \hat t$ has been replaced by $t \triangleq v_s \hat  t$ such that the sound speed in \eqref{eq:wave} becomes one.
In the case of a circular measurement geometry one assumes that the initial pressure $\source$
vanishes outside the disc $D_R\triangleq \{\rr\in\R^2 \mid \norm{\rr}<R\}$, and  the measurement sensors are located on the boundary $\partial D_R$.
In a complete data situation, the PAT image reconstruction problem consists in the recovery of the function $\source$ from the data
\begin{equation}\label{eq:data}
    (\wave \source)(\rrs,t) \triangleq p(\rrs,t)
    \quad \text{for } (\rrs,t)\in\partial D_R \times [0,T] \,,
\end{equation}
where $T>0$ is the final measurement time and $p$ denotes the solution of \eqref{eq:wave}. In practical applications, the acoustic data $\wave  \source$
are only known for a finite number of detector locations $\rrs_1, \dots , \rrs_M \in \curve$ on the measurement curve
$\curve \subseteq \partial D_R$. Additionally, one faces with practical  issues such as finite bandwidth of the
detection system, acoustic attenuation and acoustic heterogeneities.
In the paper  we  allow  a small number of detector locations (sparse data issue)  on a possible non-closed
measurement curve $\curve \subseteq \partial D_R$ (limited view issue) .

For complete measurement data  of the form \eqref{eq:data},  several efficient methods to recover $\source$ exists; see for example \cite{finch2004determining,Hal13a,jaeger2007fourier,Kun07a,xu2005universal}{}.
In the present work, we  use the filtered back projection (FBP) formula derived in \cite{FinHalRak07} for the
reconstruction of the initial pressure, which reads
\begin{equation}\label{eq:FBP}
    \forall \rr \in D_R \colon
    \quad
    \source(\rr) = \fbp (\wave \source) (\rr) \triangleq -\frac{1}{\pi R}\int_{\partial D_R}
    \Int{\abs{\rr-\rrs}}{\infty}{\frac{(\partial_t
    \wave \source)(\rrs,t)}{\sqrt{t^2-\abs{\rr-\rrs}^2}}}{t}\mathrm{d}S(z).
\end{equation}
Note that  \eqref{eq:FBP} assumes  data for all $t>0$. In  the numerical results we truncate the inner integration in \eqref{eq:FBP}
at  final measurement time $t_{\rm end}$ such that all singularities of the initial pressure have passed the measurement
curve until $t_{\rm end}$.
Additionally, \eqref{eq:FBP}  is  discretized in $\rr$, $\rrs$ and $t$ and the resulting discretization of $\fbp$ will be denoted by $\FBP \colon \R^{M \times N} \to \R^{d \times s} $.
While in the time variable we  discretize sufficiently fine, the number $M$
of sensor locations  will be small, resulting in a so-called sparse data problem.

Since we need a separate sensor for each spatial measurement sample, the number of spatial samples $\rrs_1, \dots ,\rrs_M \in \partial D_R$ is limited. Recently, systems with 64 line detectors have been demonstrated to work~\cite{gratt201564line,paltauf2017piezoelectric}{}. To keep costs low, the number of detectors will still be limited in the future and smaller than required
for artifact-free imaging advised by Shannon's sampling theory \cite{haltmeier2016sampling}{}.
This results in highly under-sampled data, which causes stripe-like artifacts in the FBP reconstruction. Our data is assumed to be sufficiently sampled
in the time domain (according Shannon's sampling theory) which is justified since time samples can easily be acquired at high sampling rate.
The   goal is to eliminate (or at least  significantly reduce) the under-sampling artifacts caused by the small number of detectors and the limited view and to improve the overall reconstruction quality. To achieve these goals we use deep learning and in particular deep CNNs.

\subsection{Deep learning image reconstruction}
\label{sec:deep}

Consider  the  following  general image reconstruction problem
\begin{equation}\label{eq:inv}
\text{Find  image } \; \Xout  \; \text{ from data } \quad
\Xin  =   \waved  (\Xout ) + \noise \,.
\end{equation}
Here $\Xout \in \R^{d \times d}$ is the image to be reconstructed,
$\Xin \in \R^{M \times N}$ are the given data,
$\waved \colon \R^{d \times d} \to \R^{M \times N} $ is the forward operator or imaging
operator, and $\noise$ models the noise.
As we show below, the image reconstruction task \eqref{eq:inv} can be seen
as a supervised machine learning problem, which can  be solved by deep neural networks.
In that context, one aims at finding  a function $\NN \colon \R^{\Nt \times \Nx}\rightarrow\R^{\Ny \times \Ny}$
that maps the  input image $\Xin  \in \R^{\Nt \times \Nx}$ to
an output image $\Xout  \in \R^{\Ny \times \Ny}$. In deep  learning,
$\NN$ is  taken as deep CNN.

For image reconstruction with deep learning, the network function
$\NN\colon \R^{\Nt \times \Nx}\rightarrow\R^{\Ny \times \Ny}$  takes the  form
\begin{equation} \label{eq:nn}
	 \NN \triangleq  \CNN  \circ \FBP \triangleq  (\sigma_L  \circ \W_L)  \circ \cdots \circ (\sigma_1  \circ \W_1) \circ \FBP  \,.
\end{equation}
Here $\FBP \colon \R^{\Nt \times \Nx}\rightarrow\R^{\Ny \times \Ny}$  maps the given data to an intermediate reconstruction in the imaging   domain. It may be taken as the adjoint of the
forward operator $\waved$; however also other choices are reasonable.
The  composition  $  \CNN =  (\sigma_L  \circ \W_L)  \circ \cdots \circ (\sigma_1  \circ \W_1) \colon \R^{\Ny \times \Ny} \to \R^{\Ny \times \Ny}$
is a layered CNN, where each
component  $\sigma_\ell  \circ \W_\ell$ is the product of a linear affine transformation (represented as a matrix) $ \W_\ell\in
\R^{D_{\ell+1}\times D_\ell}$ (where we leave the affine term $b_\ell\in\R^{D_{\ell+1}}$ out for better readability)
and a nonlinear function $\sigma_\ell  \colon \R \to \R$  that is applied component-wise.
Here $L$ denotes the number of layers, $\sigma_\ell$
are so called activation functions and $\WW \triangleq (\W_1, \ldots,\W_L)$ is the weight vector.
In CNNs, the weight matrices $\W_\ell$ are block diagonal, where each block corresponds to  a convolution  with a filter of small support and the number of blocks corresponds to the number of different filters (or channels) used in each layer. Each block is therefore a sparse band matrix,  where the non-zero entries of the band matrices determine the filters of the convolution and the number of different diagonal bands corresponds to the number of channels in the previous layer.
Modern neural networks also use additional types of operations, for example max-pooling and concatenation layers, and more  complex network architectures~\cite{goodfellow2016deep}{}.
The image reconstruction network  \eqref{eq:nn} can easily be extended to such CNNs using  more general structures.

In order to adjust the reconstruction function $\NN$  to a particular reconstruction problem and  phantom class,  the weight vector   $\WW$ is selected depending on a set of training data $\TT \triangleq \{(\Xin_n,\Xout_n)\}_{n=1}^{\Ntrain}$.
For this  purpose,  the weights are adjusted in such a way, that the overall  error  of $\NN$ made on the training set is small.
This is achieved by minimizing the error function
\begin{equation} \label{eq:err}
 	E(\TT; \NN) \triangleq
 	\sum_{n=1}^{\Ntrain} d(\NN(\Xin_n),\Xout_n) \,,
\end{equation}
where $d\colon \R^{\Nt \times \Nx} \times \R^{\Nt \times \Nx} \to [0,\infty)$ is a
distance measure  that quantifies the error made  by the network function on the $n$th training sample
$(\Xin_n,\Xout_n)$. Typical choices for the used  distance measure   (or loss function)
are the mean absolute or the mean squared metric.
During the training phase, the weights $\WW$ are adjusted such that the error $E$ is minimized.
Standard methods for that purpose are based on stochastic gradient descent.

\subsection{Proposed reconstruction networks}
\label{sec:dpat}

For PAT, in the general image reconstruction problem \eqref{eq:inv}, the data $\Xin$ are
the discrete measured PAT data  and the output  $\Xout$ is the
discretized initial  pressure $\source$ in \eqref{eq:wave}.  The forward  problem
$\waved \colon \R^{d \times d} \to \R^{M \times N} $  is the discretized solution  operator of the wave
equation, evaluated at a small number of spatial detector locations.
To recover $\Xout$ from $\Xin$ we will use a reconstruction network of the
form \eqref{eq:nn}, where    $\FBP$ is the  discretization of the FBP
operator $\fbp$. The reconstruction network   \eqref{eq:nn} then can be interpreted
to first calculate an intermediate reconstruction by  applying the FBP algorithm
to the  data $\Xin$. The intermediate reconstruction contains under-sampling and limited view artifacts
that are removed by the subsequent  neural network function $\CNN$ applied in the second step.

There are many specialised CNN designs for various tasks.  In this paper we use two different
networks namely, the U-Net \cite{ronneberger2015unet} and, for comparison purpose,
a simple CNN that we name S-Net.
\begin{itemize}
\item \textbf{U-Net:}
In this case the reconstruction network \eqref{eq:nn} takes the form $\NN = \UNET \circ \FBP $, where  $\UNET $ is the
U-Net. The U-Net was initially proposed for image segmentation in \cite{ronneberger2015unet}{}, and lately has been used successfully for reconstruction tasks like low dose CT~\cite{han2016deep,jin2017deep}
and PAT \cite{antholzer2017deep}{}.  The U-Net is a deep CNN, where each convolution
is followed by the same nonlinearity, namely the rectified linear unit (ReLU) which is defined by $\ReLU(x)\triangleq \max\{x,0\}$.
Since the structure of the residual image $\Xout - \FBP(\Xin)$ is often simpler than the structure of the original image,
we employ residual learning \cite{han2016deep}{}.
This means that we add $\FBP (\Xin)$ to the output of our NNs. Thus the error gets small if the
output of the NN is close to $\Xout -\FBP (\Xin)$.

\item \textbf{S-Net:}
The simple  network that we use for comparison purposes is based on a layered CNN
that only consists  of  three  layers. The proposed   S-Net takes the form
\begin{equation}\label{eq:SimpleNN}
    \SNET = \W_3 \circ \sigma \circ \W_2 \circ \sigma \circ \W_1 \circ \FBP \,,
\end{equation}
where all  convolutions $\W_\ell \colon \R^{\Ny \times \Ny \times D_\ell  } \to \colon \R^{\Ny \times \Ny \times D_{\ell+1} }  $
are selected to have a kernel of size $(7,7)$, and the number $D_\ell$ of channels is 64, 32 and 1 for
the first, second and third layer, respectively. The nonlinearity is taken as $\sigma=\ReLU$ for all layers
and we use zero-padding before each convolution in order to have the same
image size in all layers.
\end{itemize}

Numerical results with the  proposed reconstruction networks are
presented in the following
section. Both network designs  yield good results. Due to its simplicity, for the S-Net this
might  be slightly surprising.

\section{Results}
\label{sec:results}

In this section we demonstrate that the deep learning  framework works well for
image reconstruction in PAT using sparse limited view data. For that purpose, we  simulate sparsely sampled PAT data  and use them for training and testing the two reconstruction networks. The results are compared with plain FBP
reconstruction  and total variation (TV) minimization.

\begin{psfrags}
\psfrag{B}{$D_R\setminus R$}
\psfrag{D}{$R$}
\psfrag{C}{$\curve$}
\begin{figure}[htb!]
\floatbox[{\capbeside\thisfloatsetup{capbesideposition={right,center},capbesidewidth=0.4\textwidth}}]{figure}[\FBwidth]
{\caption{\textbf{Measurement geometry used for the  numerical experiments.} The acoustic pressure is observed at 24 sensor positions
$\rrs_1, \dots, \rrs_{24}$  (indicated by the green dots)  that are located  on the non-closed measurement curve
$\curve = \set{\rrs \colon \norm{\rrs}_2 = \SI{50}{mm} \wedge \rrs_2 < \SI{11.1}{mm}}$
forming a circular arc. The phantoms to be reconstructed are contained in the rectangular domain
$[-\SI{10}{mm} ,\SI{10}{mm} ] \times [-\SI{20}{mm} ,\SI{5}{mm} ]$ that has parts 
outside the stability region $R$
(defined as the convex hull of the measurement curve). The corresponding PAT image reconstruction problem
is the combination of a sparse data (small number of sensors) and a limited view (non-closed  measurement curve)
problem.  \label{fig:setup}}}
{ \includegraphics[width=0.5\textwidth]{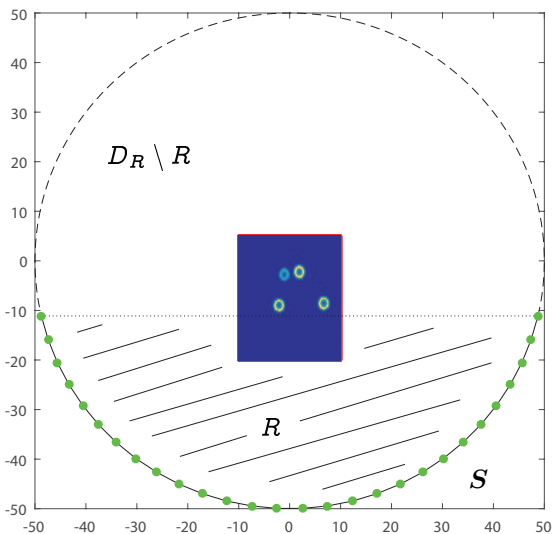}\hspace{0.02\textwidth}}
\end{figure}
\end{psfrags}

\subsection{Training and  evaluation data}

\begin{figure}[htb]
    \center
    \includegraphics[width=\textwidth]{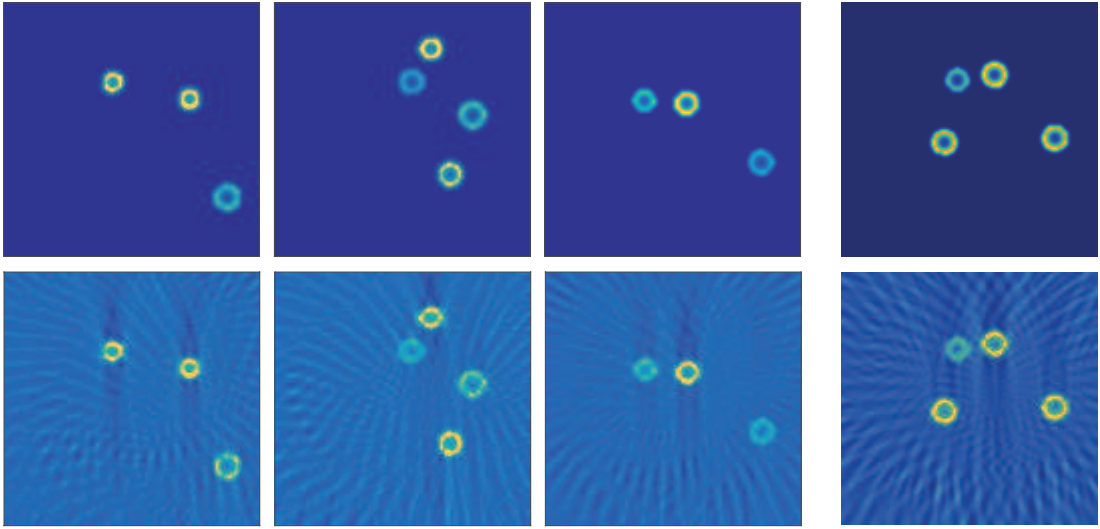}
    \caption{\textbf{Training and evaluation data.} Examples of randomly generated combinations of ring-shaped phantoms
    (top)  and  corresponding FBP reconstructions (bottom) containing under-sampling artifacts. The left three images contain examples from the training set; the
    right image  is used  for evaluation and is not part of the training data.
    In the FBP reconstructions  one clearly sees the typical under-sampling artifacts.}
    \label{fig:data}
\end{figure}

For the  presented result we generated $\Ntrain = 1000$ training phantom images $\Xout_n \in \R^{128 \times 128}$ and
$200$  evaluation phantom images $\Xout_{n'} \in \R^{128 \times 128}$, each being
a  discretization of some initial pressure $\source$ in \eqref{eq:wave} of size $128\times128$ on the domain
$[-\SI{10}{mm} ,\SI{10}{mm} ] \times [-\SI{20}{mm} ,\SI{5}{mm} ]$.
These phantoms were created by randomly superimposing  2 to 6 non-negative ring-shaped phantoms with random positions, sizes and magnitudes.
We normalize each phantom image such that it has maximal intensity value equals one.
See  Figure \ref{fig:data} for some  examples, where the  three images on left hand side are from the training  set, and
the right hand side image  is used for  evaluation. The bottom images in Figure \ref{fig:data} show the corresponding reconstructions with
the FBP algorithm, in which one can clearly see under-sampling artifact.
All radial profiles are smooth and contain blur modeling  the point spread function of the PAT imaging system and exponential decay modeling   the decrease of optical energy  within the light absorbing structures.
In particular, such structures allow to investigate the performance   of the proposed deep learning
methods on phantom classes without sharp boundaries. Results for  piecewise constants phantoms  can be found in \cite{antholzer2017deep}{}.

For the numerically  simulated acoustic data $\Xin = \waved (\Xout)$ we use
$M = 28$ uniformly  spaced detector locations arranged along
the measurement  curve
\begin{equation*}
    \curve = \set{\rrs \colon \norm{\rrs}_2 = \SI{50}{mm} \wedge \rrs_2 < \SI{-11.1}{mm}} \,.
\end{equation*}
The arrangement of the sensors along $\curve$ and the relative location of the imaging domain
are shown  in  Figure~\ref{fig:setup}.  Such a setup combines the sparse data  problem with the
limited view problem.  We take  $N = 2963$  time samples taken uniformly in the interval
$[0,  \SI{67.3}{mm}]$ (the re-scaled final measurement $\hat t_{\rm end} v_s = \SI{67.3}{mm}$  corresponds  to a physical time
of $\hat{t}_{\rm end} =\SI{44.9}{\mu s}$). To demonstrate   stability  of our deep learning approach with respect to measurement
error, we added \SI{10}{\percent} Gaussian  white noise to the simulated PAT measurements.

\begin{figure}[htb]
    \center
    \includegraphics[width=\textwidth]{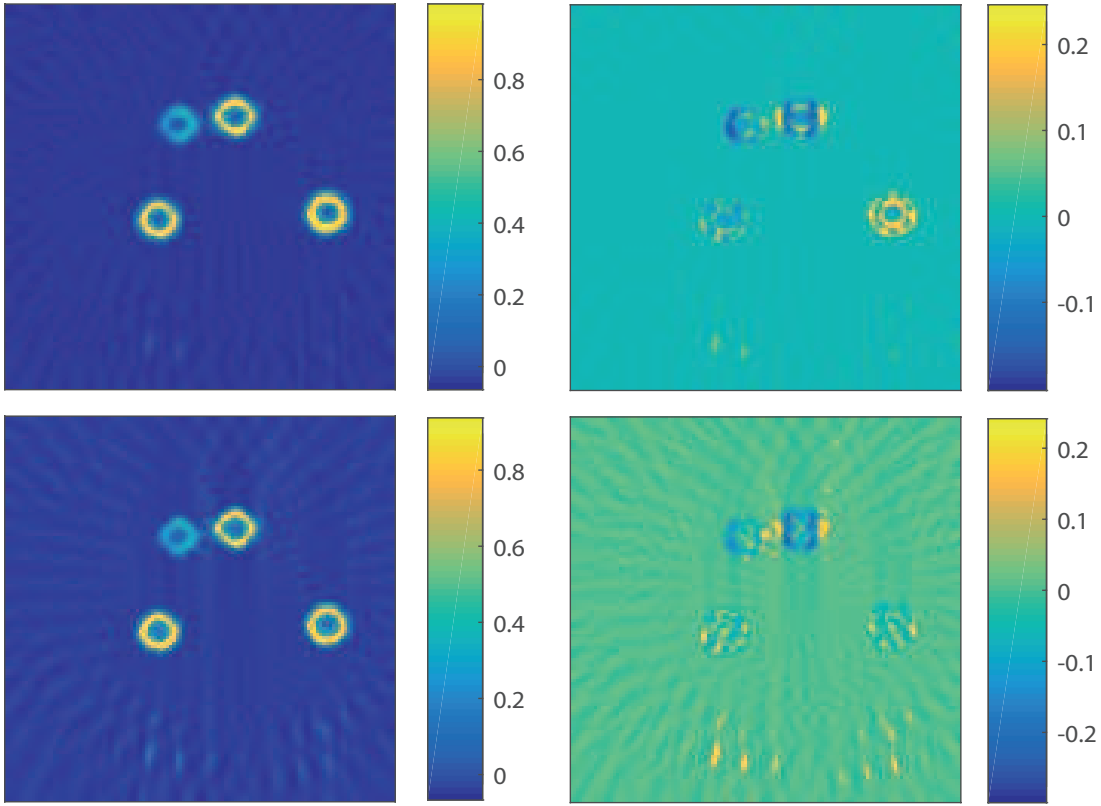}
    \caption{\textbf{Results with the reconstruction networks.}
    Top: Reconstructed image using the U-Net (left) and difference to the true phantom (right); the relative mean squared error is 0.026. Bottom:
    Reconstructed phantom using the S-net (left) and difference to the true phantom
     (right); the relative mean squared error is 0.33.}
    \label{fig:results}
\end{figure}

\subsection{Network optimization}

In order to optimize the  reconstruction networks on the training data set $ \TT = \{(\Xin_n,\Xout_n)\}_{n=1}^{\Ntrain}$, we use the mean absolute error as  error metric in \eqref{eq:err}. This   yields to  the problem of finding the parameter vector $\WW$  of the network
function $\NN \triangleq \NN_{\WW}$  by minimizing
\begin{equation} \label{eq:err2}
 	E(\TT; \NN_{\WW}) \triangleq
 	\sum_{n=1}^{\Ntrain} \norm{\Xout_n - \NN_{\WW}(\Xin_n) }_1    \,.
\end{equation}
For that purpose we did not use standard gradient  descent, because
evaluation of the full gradient is time consuming. Instead, we use stochastic batch gradient descent for minimizing
 \eqref{eq:err2}. In batch stochastic gradient descent, at each sweep (a cycle of iterations),
one partitions the training  data  into small random  subsets of equal size (batch size)  and then performs a gradient step with respect to each subset. To be precise, the update rule is given by
\begin{equation}\label{eq:SGD}
    \WW^{(k+1)} = \WW^{(k)} - \eta \nabla_{\WW} E(\TT^{(k)}, \NN_{\WW^{(k)}}) + \beta (\WW^{(k)}-\WW^{(k-1)}) \,,
\end{equation}
where $\TT^{(k)}\subseteq \TT$ is the training batch of the $k$th iteration and $\eta$ is called learning  rate and $\beta$ momentum.
Still, the training   procedure can be costly. We emphasize, however,  that the optimization
has to be performed only once and prior to the actual image reconstruction.
After training, the  weights $\WW$ are fixed and can be used to evaluate the reconstruction network
 $\NN  = \NN_{\WW}$   on new PAT data.
In the present study, to optimize our networks we use a batch size of 1 and take
$\eta = 0.001$ and $\beta=0.99$.

\begin{figure}[htb]
    \center
    \includegraphics[width=\textwidth]{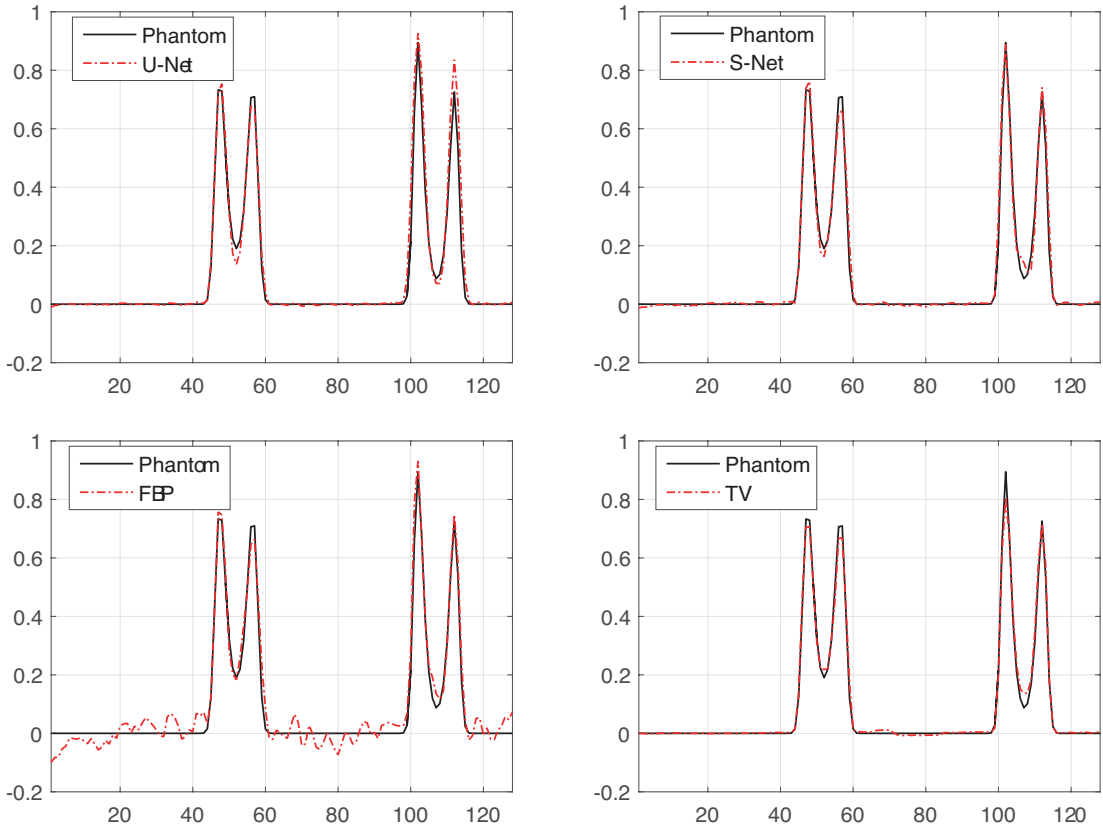}
    \caption{\textbf{Horizontal cross sections.} The images show horizontal cross section  through  the upper two rings
    comparing the original phantom  with U-Net reconstruction (top left), the S-Net reconstruction  (top right), the FBP reconstruction  (bottom left) and
    the TV-minimization (bottom right).}
    \label{fig:crosssection}
\end{figure}

\begin{figure}[htb]
    \center
    \includegraphics[width=\textwidth]{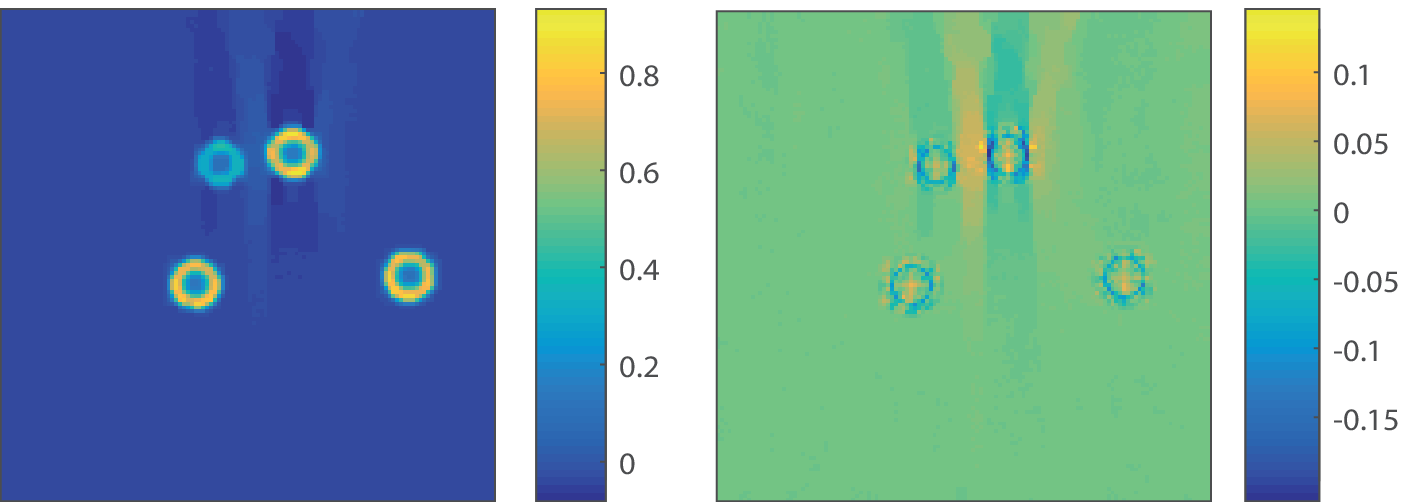}
     \caption{\textbf{Results with TV-minimization.}
     Reconstruction using TV-minimization (left) and difference to the true phantom (right).
     The relative mean squared error is 0.016.}
    \label{fig:TV}
    \end{figure}

\subsection{Reconstruction  results}

Results  with the  reconstruction networks (the U-net and the S-net)
are shown in  Figure~\ref{fig:results}.  We see that both networks are able to remove most of the  under-sampling artifacts.
The more complex U-Net  gives  better results, but also takes longer to train and apply. Taking a look at a horizontal cross section of
the reconstructed phantom (Figure~\ref{fig:crosssection}) we can see that both networks overestimate the minimal values within the ring structures.
This suggests that classes of highly oscillating phantoms might be quite challenging to reconstruct.
Further research is required to find out how to handle such cases.
It is still surprising that the second NN, which is quite simple, performs that
well. However, we expect that the S-net might  struggle with more complex phantom classes. Anyway these results  encourage
to design  new and well-suited
networks for PAT image  reconstruction.

For comparison purpose, we also tested TV minimization for image reconstruction,
\begin{equation}  \label{eq:TV}
	 \frac{1}{2} \norm{ \waved (\Xin) -  \Xout}_2^2
	+ \lambda \, \sum_{\mathbf{i}}  \sqrt{ \abs{ (\Dnum_1 \Xout)_{\mathbf{i}} }^2 + \abs{ (\Dnum_2 \Xout)_{\mathbf{i}} }^2 }    \to \min_{\Xout}\,,
\end{equation}
where $\Dnum = [\Dnum_1, \Dnum_2]$ is the discrete gradient operator
and $\lambda$ the regularization  parameter.   For solving \eqref{eq:TV}
we use the  algorithm proposed in  \cite{sidky2012convex}
(an instance of the Pock-Chambolle algorithm \cite{chambolle2011first}
for TV minimization).
Since the considered phantoms do not contain structures at different scales,
TV-minimization method performs  well, and in fact
shows the lowest error relative mean squared error
$\MSE(\Xout)  \triangleq \norm{\Xout- \Xout_{\rm rec}}^2 /\norm{\Xout}_2^2$, where $\norm{\Xout}_2$ denote the $\ell^2$-norm of $\Xout$.
However, TV-minimization requires choosing a good regularization  parameter and performing a relatively large number of iterations.
For the  results shown in Figure \ref{fig:TV} we have chosen the optimization parameter $\lambda=0.005$ by hand and, in order to  get small $\ell^2$-reconstruction error,  performed  50 iterations. We note that the TV-minimization introduces new additional artifacts around the pair of rings which are relatively close together.

\subsection{Computational resources}

We used Keras~\cite{keras} with TensorFlow~\cite{tensorflow}
to train an evaluate the proposed reconstruction networks (U-Net and S-Net). The FBP and the TV-minimization
are implemented in MATLAB.
We ran all our experiments on a computer using an Intel i7-6850K and an NVIDIA 1080Ti.
To iterate through our entire training set we need \SI{5}{s} for the S-Net~\eqref{eq:SimpleNN} and  \SI{16}{s}
for the U-Net. Evaluating 100 sample images requires \SI{0.9}{s}  for the S-Net  and  \SI{3.1}{s}   for the U-Net.
Hence  a single image is reconstructed in a fraction of seconds by both methods; for example the S-net reconstructs
the $128\times 128$  images at  \SI{111}{\hertz} rate.  One application with the current implementation of the discrete FBP  takes \SI{0.35}{s}
and the used TV-minimization  algorithm  needs \SI{0.94}{s} for one iteration. This results in an overall
reconstruction time for TV minimization around \SI{45}{s}.
Since the latter algorithms are implemented in MATLAB and do not use the GPU,
the comparison of computation times is not completely fair and there is  room for accelerating the  FBP algorithm and TV-minimization.
Especially, the  FBP algorithm in  combination with CNNs, both implemented  on
GPUs, will  give high resolution artifact-free reconstructions in real time.

\section{Conclusion}
\label{sec:conclusion}

In this paper we proposed  a deep learning framework for image reconstruction in PAT using sparse data including the  limited
view setting.  The proposed reconstruction  structure \eqref{eq:nn} consists in  first applying the FBP algorithm and then using a CNN to remove artifacts.
For the used CNN we investigated the established U-Net as well as the  simple S-Net. Both of the proposed networks are able to improve  the overall image quality. As expected, the more complex U-Net yields better result and offers a reconstruction quality  comparable to iterative TV-minimization.
Both reconstruction networks can be applied in real time, in contrast to iterative methods which are much slower.
In future research, we will use more complex simulated and real-world PAT data, where we expect our method to be faster  and
comparable  to TV-minimization in terms of reconstruction error. We also think it is beneficial to include adjustable weights in the backprojection step,  which in the moment  just used the FBP algorithms with
prescribed weights.

\section*{Acknowledgement}

SA and MH  acknowledge support of the Austrian Science Fund (FWF), project P 30747. The work of RN has been supported by the FWF, project P 28032.

\end{document}